%% file: _gogmagog.tex
\documentclass[a4paper]{article}
\input{_gogmagog_head}

\begin{document}
\maketitle

\begin{abstract}
A sub-problem of the open problem of finding an explicit bijection between alternating sign matrices and totally symmetric self-complementary plane partitions consists in finding an explicit bijection between so-called $(n,k)$ Gog trapezoids and $(n,k)$ Magog trapezoids. A quite involved bijection was found by Biane and Cheballah in the case $k=2$. We give here a simpler bijection for this case. 

\bigskip

\noindent\textbf{Key words and phrases:} bijection, Gog, Magog, alternating sign matrix, totally symmetric self-complementary plane partition.

\medskip
\noindent\textbf{AMS classification:} 05A19.

\end{abstract}


\input{_gogmagog_body}

\bibliographystyle{alpha}
\bibliography{main}
\end{document}

%% file: _gogmagog_head.tex
\usepackage[textwidth=15cm,textheight=22cm]{geometry}
\usepackage{amssymb}
\usepackage{amsmath,cases}
\usepackage{amsthm}
\usepackage{enumitem}
\usepackage{xcolor}
\usepackage{graphicx}
\usepackage{textcomp}
\usepackage{mathrsfs}
\usepackage{bbm}
\usepackage{bm}
\usepackage{float}
\usepackage{caption}
\usepackage{chemarrow}
\usepackage{stmaryrd}
\usepackage{url}
\usepackage{calc}
\usepackage{ifthen}
\usepackage{psfrag}
\usepackage[active]{srcltx}   
\usepackage[all,cmtip]{xy}
\usepackage{pdfsync}																				
\usepackage{yhmath}
\usepackage[colorlinks=true,bookmarksnumbered=true, pdfauthor={Jérémie Bettinelli}]{hyperref}									

\definecolor{violet}{RGB}{255,85,255}
\definecolor{vert}{RGB}{142,255,142}
\definecolor{jaune}{RGB}{255,255,113}
\definecolor{rouge}{RGB}{255,56,56}
\definecolor{bleu}{RGB}{142,142,255}
\newcommand{\col}[2]{\colorbox{#1}{\vphantom1\ensuremath{#2}}}

\newcommand{\hr}{\texorpdfstring}

\newcommand{\de}{\mathrel{\mathop:}\hspace*{-.6pt}=}

\newcommand{\calG}{\mathcal{G}}

\newcommand{\M}{\mathcal{M}}

\definecolor{gris}{gray}{0.7}
\definecolor{grisf}{gray}{0.4}

\theoremstyle{plain}
\newtheorem{thm}{Theorem}

\newtheorem{defi}{Definition}

\theoremstyle{definition}

\newtheorem*{ack}{Acknowledgment}

\newenvironment{pre}[1][\proofname]{%
  \proof[#1]%
}{\endproof}

\makeatletter
\renewcommand*{\@fnsymbol}[1]{\ensuremath{\ifcase#1\or \dagger\or \ddagger\or
   \mathsection\or \mathparagraph\or \|\or **\or \dagger\dagger
   \or \ddagger\ddagger \else\@ctrerr\fi}}
\makeatother
\title{A simple explicit bijection between $(n,2)$-Gog and Magog trapezoids}
\author{J\'er\'emie \textsc{Bettinelli}\thanks{CNRS \& Laboratoire d'Informatique de l'\'Ecole polytechnique; \href{mailto:jeremie.bettinelli@normalesup.org}{\nolinkurl{jeremie.bettinelli@normalesup.org}}; \nolinkurl{www.normalesup.org/}\texttildelow\nolinkurl{bettinel}.}}

%% file: _gogmagog_body.tex
\section{Introduction}

One of the most famous open problems in bijective combinatorics is to find an explicit bijection between alternating sign matrices of a given size and totally symmetric self-complementary plane partitions of the same size. These objects of combinatorial interest have been known since the end of the '90s to be equinumerous \cite{Andrews94,Zeilberger96} but, as of today, there is no direct bijective proof of this fact. We refer the reader to~\cite{bressoud,cheballah11these} for more information on this story.

The previous objects are in known bijections with arrays of integers called \emph{Gog} and \emph{Magog} triangles. These triangles are Young diagrams of shape $(n,{n-1},\dots,\allowbreak 2,1)$ for some positive integer $n$, and they are filled with positive integers satisfying monotonicity conditions along vertical, horizontal and possibly diagonal lines. Although they satisfy very similar monotonicity conditions, nobody managed to find a direct bijection between these integer-filled triangles so far. Another surprising fact is that, if we only consider the first $k$ rows of the triangles, the objects we obtain are also equinumerous. These objects, called \emph{$(n,k)$ trapezoids}, were introduced in~\cite{MillsRobbinsRumsey}, where they were conjectured to be equinumerous. This was later proved by Zeilberger \cite[Lemma~1]{Zeilberger96}.

The supposedly simplest problem of finding an explicit bijection between $(n,k)$-Gog trapezoids and $(n,k)$-Magog trapezoids has been solved only for $k\le 2$. In fact, for $k=1$, the objects are exactly the same, so there is nothing to prove. There is, however, a refined conjecture appearing in~\cite{MillsRobbinsRumsey}, which involves the number of entries equal to~$1$ and the number of entries equal to the maximum possible value in the first and last rows. Proving this conjecture in a bijective way, even in the case $k=1$ is nontrivial; it was achieved by Krattenthaler~\cite{krattenthalerGogMagog}. In fact, Krattenthaler conjectured that more general objects called $(\ell,n,k)$ trapezoids are equinumerous and that the previous statistics coincide; his bijection was established in this more general context, that is, between $(\ell,n,1)$-Gog and Magog trapezoids.

For $k=2$, a bijection was found by Biane and Cheballah~\cite{bianecheballah12}. Their bijection is relatively complicated and uses the so-called Sch\"utzenberger involution. It does not match the previous statistics of Mills, Robbins, and Rumsey. It does, however, match different statistics, expressed in terms of the rightmost entry for a Gog trapezoid and in terms of the two rightmost entries of both rows for Magog trapezoids. 
In this note, we give a different bijection for this case. Our bijection is very simple and involves only one operation. It does not match either aforementioned statistics.

\begin{ack}
I thank Jean-Fran\c{c}ois Marckert for introducing this problem to me.
\end{ack}

\section{Magog and Gog trapezoids}

\begin{figure}[ht]
		\psfrag{1}[][]{$1$}
		\psfrag{2}[][]{$2$}
		\psfrag{3}[][]{$3$}
		\psfrag{4}[][]{$4$}
		\psfrag{5}[][]{$5$}
		\psfrag{6}[][]{$6$}
		\psfrag{7}[][]{$7$}
		\psfrag{8}[][]{$8$}
		
		\psfrag{q}[B][B][.8]{\textcolor{red}{$g_{1,1}$}}
		\psfrag{s}[B][B][.8]{\textcolor{red}{$g_{1,2}$}}
		\psfrag{d}[B][B][.8]{\textcolor{red}{$g_{1,3}$}}
		\psfrag{f}[B][B][.8]{\textcolor{red}{$g_{1,4}$}}
		\psfrag{g}[B][B][.8]{\textcolor{red}{$g_{1,5}$}}
		\psfrag{h}[B][B][.8]{\textcolor{red}{$g_{1,6}$}}
		\psfrag{j}[B][B][.8]{\textcolor{red}{$g_{1,7}$}}
		\psfrag{k}[B][B][.8]{\textcolor{red}{$g_{1,8}$}}
		
		\psfrag{a}[B][B][.8]{\textcolor{red}{$g_{2,1}$}}
		\psfrag{z}[B][B][.8]{\textcolor{red}{$g_{2,2}$}}
		\psfrag{e}[B][B][.8]{\textcolor{red}{$g_{2,3}$}}
		\psfrag{r}[B][B][.8]{\textcolor{red}{$g_{2,4}$}}
		\psfrag{t}[B][B][.8]{\textcolor{red}{$g_{2,5}$}}
		\psfrag{y}[B][B][.8]{\textcolor{red}{$g_{2,6}$}}
		\psfrag{u}[B][B][.8]{\textcolor{red}{$g_{2,7}$}}
		
		\psfrag{Q}[B][B][.8]{\textcolor{red}{$m_{1,1}$}}
		\psfrag{S}[B][B][.8]{\textcolor{red}{$m_{1,2}$}}
		\psfrag{D}[B][B][.8]{\textcolor{red}{$m_{1,3}$}}
		\psfrag{F}[B][B][.8]{\textcolor{red}{$m_{1,4}$}}
		\psfrag{G}[B][B][.8]{\textcolor{red}{$m_{1,5}$}}
		\psfrag{H}[B][B][.8]{\textcolor{red}{$m_{1,6}$}}
		\psfrag{J}[B][B][.8]{\textcolor{red}{$m_{1,7}$}}
		
		\psfrag{A}[B][B][.8]{\textcolor{red}{$m_{2,1}$}}
		\psfrag{Z}[B][B][.8]{\textcolor{red}{$m_{2,2}$}}
		\psfrag{E}[B][B][.8]{\textcolor{red}{$m_{2,3}$}}
		\psfrag{R}[B][B][.8]{\textcolor{red}{$m_{2,4}$}}
		\psfrag{T}[B][B][.8]{\textcolor{red}{$m_{2,5}$}}
		\psfrag{Y}[B][B][.8]{\textcolor{red}{$m_{2,6}$}}
		\psfrag{U}[B][B][.8]{\textcolor{red}{$m_{2,7}$}}
		\psfrag{I}[B][B][.8]{\textcolor{red}{$m_{2,8}$}}		
		
		\psfrag{n}[][]{$(8,2)$-Gog trapezoid}
		\psfrag{m}[][]{$(8,2)$-Magog trapezoid}
	\centering\includegraphics[width=.95\linewidth]{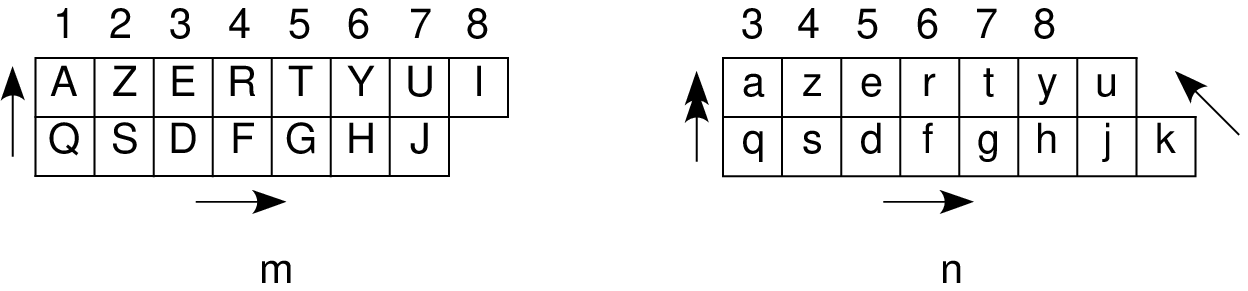}

		\psfrag{a}[B][B]{\textcolor{red}{$2$}}
		\psfrag{z}[B][B]{\textcolor{red}{$2$}}
		\psfrag{e}[B][B]{\textcolor{red}{$4$}}
		\psfrag{r}[B][B]{\textcolor{red}{$5$}}
		\psfrag{t}[B][B]{\textcolor{red}{$6$}}
		\psfrag{y}[B][B]{\textcolor{red}{$7$}}
		\psfrag{u}[B][B]{\textcolor{red}{$8$}}
		
		\psfrag{q}[B][B]{\textcolor{red}{$1$}}
		\psfrag{s}[B][B]{\textcolor{red}{$1$}}
		\psfrag{d}[B][B]{\textcolor{red}{$2$}}
		\psfrag{f}[B][B]{\textcolor{red}{$4$}}
		\psfrag{g}[B][B]{\textcolor{red}{$4$}}
		\psfrag{h}[B][B]{\textcolor{red}{$5$}}
		\psfrag{j}[B][B]{\textcolor{red}{$7$}}
		\psfrag{k}[B][B]{\textcolor{red}{$7$}}
		
		\psfrag{A}[B][B]{\textcolor{red}{$1$}}
		\psfrag{Z}[B][B]{\textcolor{red}{$2$}}
		\psfrag{E}[B][B]{\textcolor{red}{$2$}}
		\psfrag{R}[B][B]{\textcolor{red}{$4$}}
		\psfrag{T}[B][B]{\textcolor{red}{$4$}}
		\psfrag{Y}[B][B]{\textcolor{red}{$6$}}
		\psfrag{U}[B][B]{\textcolor{red}{$7$}}
		\psfrag{I}[B][B]{\textcolor{red}{$7$}}
		
		\psfrag{Q}[B][B]{\textcolor{red}{$1$}}
		\psfrag{S}[B][B]{\textcolor{red}{$1$}}
		\psfrag{D}[B][B]{\textcolor{red}{$2$}}
		\psfrag{F}[B][B]{\textcolor{red}{$4$}}
		\psfrag{G}[B][B]{\textcolor{red}{$4$}}
		\psfrag{H}[B][B]{\textcolor{red}{$5$}}
		\psfrag{J}[B][B]{\textcolor{red}{$7$}}
		\psfrag{K}[B][B]{\textcolor{red}{$7$}}
		
		\psfrag{n}[][]{}
		\psfrag{m}[][]{}
	\vspace*{5mm}
	\includegraphics[width=.95\linewidth]{def}
	\vspace*{-7mm}
	\caption{The conditions satisfied by $(n,2)$-Magog and Gog trapezoids. Every sequence formed by numbers obtained by following the direction of a simple-arrowhead (respectively a double-arrowhead) arrow is non-decreasing (respectively increasing).}
	\label{def}
\end{figure}

In this note, we are solely considering $(n,2)$ trapezoids, and we furthermore impose that $n\ge 3$ in order to avoid trivialities. Let us give proper definitions (see Figure~\ref{def} for more graphical definitions and examples).
\begin{defi}\label{defm}
Let $n$ be an integer $\ge 3$. 
An \emph{$(n,2)$-Magog trapezoid} is an array of $2n-1$ positive integers $m_{1,1}$, \ldots, $m_{1,n-1}$, $m_{2,1}$, \ldots, $m_{2,n}$ such that
\begin{enumerate}[label=(\textit{\roman*})]
	\item $m_{i,j}\le m_{i,j+1}$ for all $i\in\{1,2\}$ and $j\in\{1,\ldots, n+i-3\}$\,;
	\item $m_{1,j}\le m_{2,j}\le j$ for all $j\in\{1,\ldots, n-1\}$ and $m_{2,n}\le n$\,.\label{defmii}
\end{enumerate}

\end{defi}
\begin{defi}\label{defg}
Let $n$ be an integer $\ge 3$. An \emph{$(n,2)$-Gog trapezoid} is an array of $2n-1$ positive integers $g_{1,1}$, \ldots, $g_{1,n}$, $g_{2,1}$, \ldots, $g_{2,n-1}$ such that
\begin{enumerate}[label=(\textit{\roman*})]
	\item $g_{i,j}\le g_{i,j+1}$ for all $i\in\{1,2\}$ and $j\in\{1,\ldots, n-i\}$\,;
	\item $g_{1,j}< g_{2,j}<j+2$ for all $j\in\{1,\ldots, n-1\}$\,;\label{defgii}
	\item $g_{1,j+1}\le g_{2,j}$ for all $j\in\{1,\ldots, n-1\}$\,.
\end{enumerate}
\end{defi}

We denote the sets of $(n,2)$-Magog and Gog trapezoids by $\M_n$ and $\calG_n$,
respectively.

\section{From Magog to Gog}\label{secmg}

\begin{figure}[H]
		\psfrag{m}[][]{$\mapsto$}
		\psfrag{+}[][]{$+1$}
		\psfrag{-}[][]{$-2$}
		\psfrag{<}[][]{$<$}
		\psfrag{=}[][]{$=$}
		
		\psfrag{A}[B][B]{$1$}
		\psfrag{Z}[B][B]{$2$}
		\psfrag{E}[B][B]{$2$}
		\psfrag{R}[B][B]{$4$}
		\psfrag{T}[B][B]{$4$}
		\psfrag{Y}[B][B]{$6$}
		\psfrag{U}[B][B]{$7$}
		\psfrag{I}[B][B]{$7$}
		
		\psfrag{Q}[B][B]{$1$}
		\psfrag{S}[B][B]{$1$}
		\psfrag{D}[B][B]{$2$}
		\psfrag{F}[B][B]{$4$}
		\psfrag{G}[B][B]{$4$}
		\psfrag{H}[B][B]{$6$}
		\psfrag{J}[B][B]{$7$}
		
		\psfrag{a}[B][B]{$2$}
		\psfrag{z}[B][B]{$3$}
		\psfrag{e}[B][B]{$4$}
		\psfrag{r}[B][B]{$4$}
		\psfrag{t}[B][B]{$6$}
		\psfrag{y}[B][B]{$7$}
		\psfrag{u}[B][B]{$7$}
		
		\psfrag{q}[B][B]{$1$}
		\psfrag{s}[B][B]{$1$}
		\psfrag{d}[B][B]{$2$}
		\psfrag{f}[B][B]{$2$}
		\psfrag{g}[B][B]{$2$}
		\psfrag{h}[B][B]{$2$}
		\psfrag{j}[B][B]{$4$}
		\psfrag{k}[B][B]{$5$}
	\centering
	\includegraphics[width=.95\linewidth]{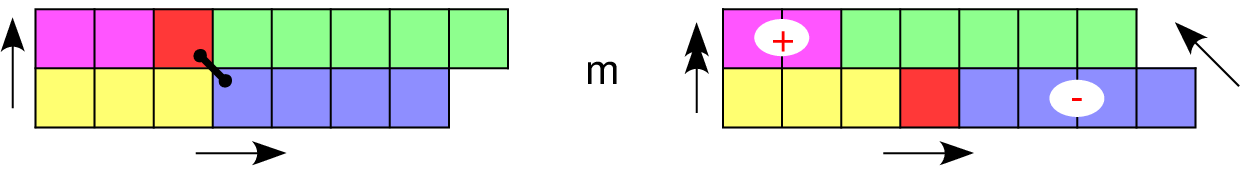}
	\vspace*{3mm}
	
	\includegraphics[width=.95\linewidth]{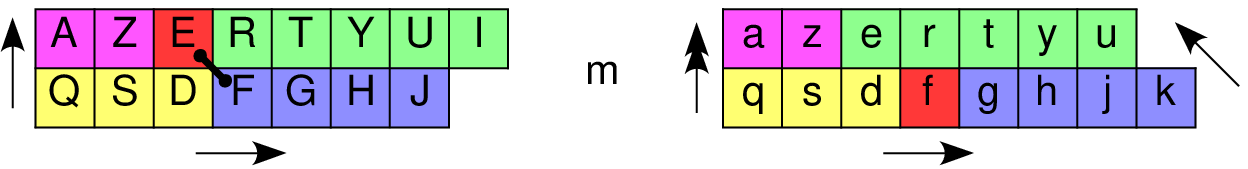}
	
		\psfrag{A}[B][B]{$1$}
		\psfrag{Z}[B][B]{$2$}
		\psfrag{E}[B][B]{$2$}
		\psfrag{R}[B][B]{$4$}
		\psfrag{T}[B][B]{$4$}
		\psfrag{Y}[B][B]{$6$}
		\psfrag{U}[B][B]{$6$}
		\psfrag{I}[B][B]{$8$}
		
		\psfrag{Q}[B][B]{$1$}
		\psfrag{S}[B][B]{$1$}
		\psfrag{D}[B][B]{$2$}
		\psfrag{F}[B][B]{$3$}
		\psfrag{G}[B][B]{$4$}
		\psfrag{H}[B][B]{$4$}
		\psfrag{J}[B][B]{$5$}
		
		\psfrag{a}[B][B]{$2$}
		\psfrag{z}[B][B]{$3$}
		\psfrag{e}[B][B]{$3$}
		\psfrag{r}[B][B]{$5$}
		\psfrag{t}[B][B]{$5$}
		\psfrag{y}[B][B]{$7$}
		\psfrag{u}[B][B]{$8$}
		
		\psfrag{q}[B][B]{$1$}
		\psfrag{s}[B][B]{$1$}
		\psfrag{d}[B][B]{$2$}
		\psfrag{f}[B][B]{$3$}
		\psfrag{g}[B][B]{$4$}
		\psfrag{h}[B][B]{$4$}
		\psfrag{j}[B][B]{$5$}
		\psfrag{k}[B][B]{$6$}
	\vspace*{8mm}
	\includegraphics[width=.95\linewidth]{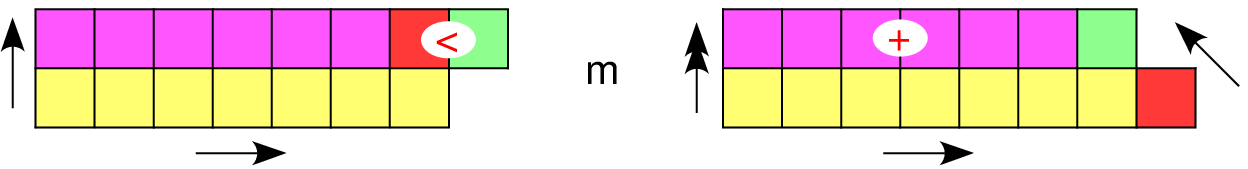}
	\vspace*{3mm}
	
	\includegraphics[width=.95\linewidth]{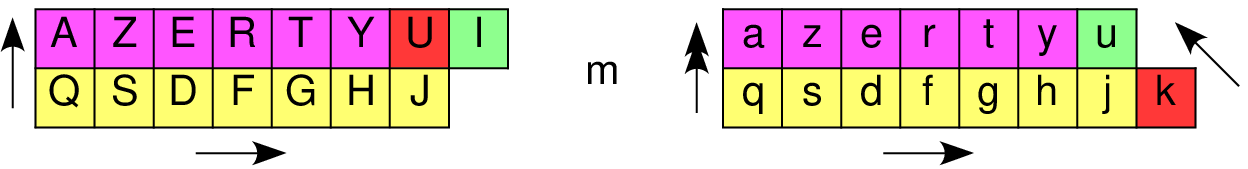}
	
		\psfrag{A}[B][B]{$1$}
		\psfrag{Z}[B][B]{$2$}
		\psfrag{E}[B][B]{$2$}
		\psfrag{R}[B][B]{$4$}
		\psfrag{T}[B][B]{$4$}
		\psfrag{Y}[B][B]{$6$}
		\psfrag{U}[B][B]{$6$}
		\psfrag{I}[B][B]{$6$}
		
		\psfrag{Q}[B][B]{$1$}
		\psfrag{S}[B][B]{$1$}
		\psfrag{D}[B][B]{$2$}
		\psfrag{F}[B][B]{$3$}
		\psfrag{G}[B][B]{$4$}
		\psfrag{H}[B][B]{$4$}
		\psfrag{J}[B][B]{$5$}
		
		\psfrag{a}[B][B]{$2$}
		\psfrag{z}[B][B]{$3$}
		\psfrag{e}[B][B]{$3$}
		\psfrag{r}[B][B]{$5$}
		\psfrag{t}[B][B]{$5$}
		\psfrag{y}[B][B]{$7$}
		\psfrag{u}[B][B]{$7$}
		
		\psfrag{q}[B][B]{$1$}
		\psfrag{s}[B][B]{$1$}
		\psfrag{d}[B][B]{$2$}
		\psfrag{f}[B][B]{$3$}
		\psfrag{g}[B][B]{$4$}
		\psfrag{h}[B][B]{$4$}
		\psfrag{j}[B][B]{$5$}
		\psfrag{k}[B][B]{$7$}
	\vspace*{8mm}
	\includegraphics[width=.95\linewidth]{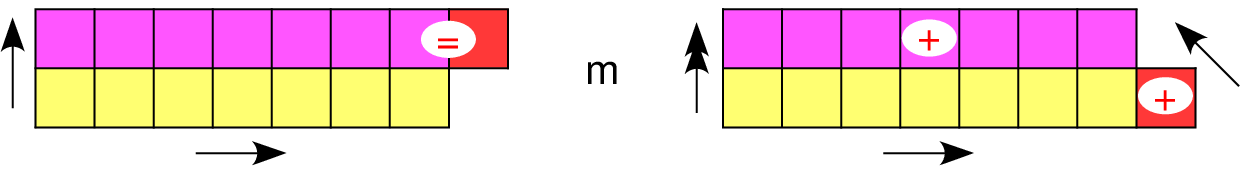}
	\vspace*{3mm}
	
	\includegraphics[width=.95\linewidth]{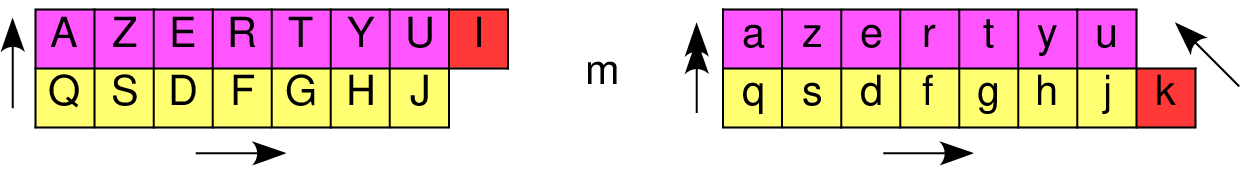}		
	\caption{The three cases of the bijection, from a Magog trapezoid to a Gog trapezoid. On the top line, the first bug is~$3$: it is symbolized by a small black diagonal line. The colored blocks are moved and, whenever there is a tag on a block, it is added to all the elements of the block.}
	\label{mag_gog}
\end{figure}

Let us consider an $(n,2)$-Magog trapezoid $M=(m_{i,j})$. We say that an integer $j\in\{1,\ldots,n-2\}$ is a \emph{bug} if $m_{1,j+1}>m_{2,j}+1$. For instance, $3$ is the only bug of the Magog trapezoid of Figure~\ref{def}. We set $\Phi_n(M)\de (g_{ij})$, where $(g_{ij})$ is constructed as follows (see Figure~\ref{mag_gog}).

\medskip
\noindent\textbf{First case: $M$ has at least one bug.}
In this case, we let $k$ be the smallest bug of~$M$ and we set
\begin{align*}
\col{violet}{g_{2,j}}&\de m_{2,j}+1\ \text{ for } 1\le j \le k-1\,; \qquad\col{vert}{g_{2,j}}\de m_{2,j+1}\ \text{ for } k\le j \le n-1\,;\\
\col{jaune}{g_{1,j}}&\de m_{1,j}\ \text{ for } 1\le j \le k\,;\qquad\col{rouge}{g_{1,k+1}}\de m_{2,k}\,; \qquad\col{bleu}{g_{1,j}}\de m_{1,j-1}-2\ \text{ for }k+2\le j \le n\,.
\end{align*}

\medskip
\noindent\textbf{Second case: $M$ does not have bugs and $m_{2,n-1}<m_{2,n}$\,.} In this case, we set
\begin{align*}
\col{violet}{g_{2,j}}&\de m_{2,j}+1\ \text{ for } 1\le j \le n-2\,; &\col{vert}{g_{2,n-1}}\de m_{2,n}\,;\\
\col{jaune}{g_{1,j}}&\de m_{1,j}\ \text{ for } 1\le j \le n-1\,; &\col{rouge}{g_{1,n}}\de m_{2,n-1}\,.
\end{align*}

\medskip
\noindent\textbf{Third case: $M$ does not have bugs and $m_{2,n-1}=m_{2,n}$\,.} In this case, we set
\begin{align*}
\col{violet}{g_{2,j}}&\de m_{2,j}+1\ \text{ for } 1\le j \le n-1\,;\\
\col{jaune}{g_{1,j}}&\de m_{1,j}\ \text{ for } 1\le j \le n-1\,; &&\col{rouge}{g_{1,n}}\de m_{2,n}+1\,.
\end{align*}

Let us check that $\Phi_n(M)\in\calG_n$. First, observe that, if~$j$ is not a bug, then, by definition, $m_{1,j+1}\le m_{2,j}+1$, so that the yellow and purple blocks always satisfy the diagonal inequalities after the mapping. It is straightforward to verify that the other inequalities are satisfied in the second and third case. In the first case, notice that $\col{jaune}{g_{1,k}}=m_{1,k}\le m_{2,k}=\col{rouge}{g_{1,k+1}}$ and $\col{rouge}{g_{1,k+1}}=m_{2,k} \le m_{1,k+1}-2 = \col{bleu}{g_{1,k+2}}$ as~$k$ is a bug. Furthermore, $\col{violet}{g_{2,k-1}}=m_{2,k-1}+1\le m_{2,k}+1 \le m_{1,k+1}-1\le m_{2,k+1}-1=\col{vert}{g_{2,k}}-1$ so that the horizontal inequalities are satisfied. Moreover, $\col{jaune}{g_{1,k}}=m_{1,k}\le m_{2,k}\le m_{1,k+1}-2 \le m_{2,k+1} -2 =\col{vert}{g_{2,k}}-2$, $\col{rouge}{g_{1,k+1}}= m_{2,k}\le m_{2,k+2}-2=\col{vert}{g_{2,k+1}}-2$, $\col{bleu}{g_{1,j}}=m_{1,j-1}-2\le m_{2,j+1}-2=\col{vert}{g_{2,j}}-2$ for $k+2\le j \le n-1$, and $\col{vert}{g_{2,j}}=m_{2,j+1}< j+2$ for $k\le j\le n-1$, so that the vertical inequalities are also satisfied. Finally, the diagonal inequalities are satisfied since $\col{rouge}{g_{1,k+1}}=m_{2,k}\le m_{2,k+1}=\col{vert}{g_{2,k}}$ and $\col{bleu}{g_{1,j}}=m_{1,j-1}-2 \le m_{2,j}-2=\col{vert}{g_{2,j-1}}-2$ for $k+2\le j \le n$.

\section{From Gog to Magog}\label{secgm}

We now consider an $(n,2)$-Gog trapezoid $G=(g_{i,j})$ and construct $\Psi_n(G)= (m_{ij})$ as follows. We define
\begin{equation}\label{defk}
k\de \max\big\{j\in\{2,\ldots,n-1\}\,: g_{2,j-1} \le g_{1,j+1}+1\big\}.
\end{equation}
This number is well defined as $g_{2,1}=2\le g_{1,3}+1$.

\begin{figure}[ht]
		\psfrag{m}[][]{$\mapsto$}
		\psfrag{+}[][]{$-1$}
		\psfrag{-}[][]{$+2$}
		\psfrag{<}[][]{$>$}
		\psfrag{=}[][]{$=$}
		
		\psfrag{A}[B][B]{$1$}
		\psfrag{Z}[B][B]{$2$}
		\psfrag{E}[B][B]{$2$}
		\psfrag{R}[B][B]{$4$}
		\psfrag{T}[B][B]{$4$}
		\psfrag{Y}[B][B]{$6$}
		\psfrag{U}[B][B]{$7$}
		\psfrag{I}[B][B]{$7$}
		
		\psfrag{Q}[B][B]{$1$}
		\psfrag{S}[B][B]{$1$}
		\psfrag{D}[B][B]{$2$}
		\psfrag{F}[B][B]{$4$}
		\psfrag{G}[B][B]{$4$}
		\psfrag{H}[B][B]{$6$}
		\psfrag{J}[B][B]{$7$}
		
		\psfrag{a}[B][B]{$2$}
		\psfrag{z}[B][B]{$3$}
		\psfrag{e}[B][B]{$4$}
		\psfrag{r}[B][B]{$4$}
		\psfrag{t}[B][B]{$6$}
		\psfrag{y}[B][B]{$7$}
		\psfrag{u}[B][B]{$7$}
		
		\psfrag{q}[B][B]{$1$}
		\psfrag{s}[B][B]{$1$}
		\psfrag{d}[B][B]{$2$}
		\psfrag{f}[B][B]{$2$}
		\psfrag{g}[B][B]{$2$}
		\psfrag{h}[B][B]{$2$}
		\psfrag{j}[B][B]{$4$}
		\psfrag{k}[B][B]{$5$}
	\centering
	\includegraphics[width=.95\linewidth]{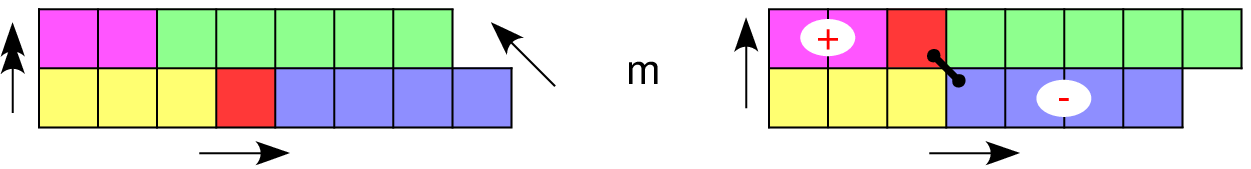}
	\vspace*{3mm}
	
	\includegraphics[width=.95\linewidth]{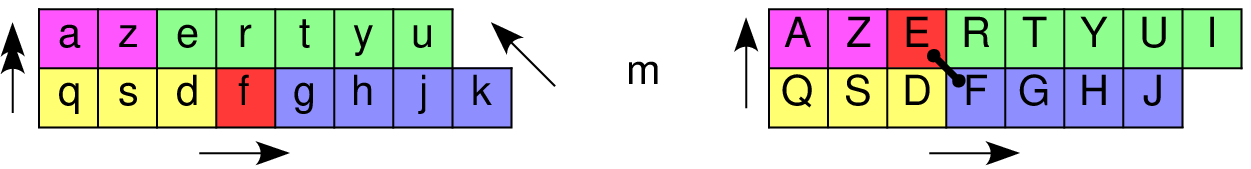}
	
		\psfrag{A}[B][B]{$1$}
		\psfrag{Z}[B][B]{$2$}
		\psfrag{E}[B][B]{$2$}
		\psfrag{R}[B][B]{$4$}
		\psfrag{T}[B][B]{$4$}
		\psfrag{Y}[B][B]{$6$}
		\psfrag{U}[B][B]{$6$}
		\psfrag{I}[B][B]{$8$}
		
		\psfrag{Q}[B][B]{$1$}
		\psfrag{S}[B][B]{$1$}
		\psfrag{D}[B][B]{$2$}
		\psfrag{F}[B][B]{$3$}
		\psfrag{G}[B][B]{$4$}
		\psfrag{H}[B][B]{$4$}
		\psfrag{J}[B][B]{$5$}
		
		\psfrag{a}[B][B]{$2$}
		\psfrag{z}[B][B]{$3$}
		\psfrag{e}[B][B]{$3$}
		\psfrag{r}[B][B]{$5$}
		\psfrag{t}[B][B]{$5$}
		\psfrag{y}[B][B]{$7$}
		\psfrag{u}[B][B]{$8$}
		
		\psfrag{q}[B][B]{$1$}
		\psfrag{s}[B][B]{$1$}
		\psfrag{d}[B][B]{$2$}
		\psfrag{f}[B][B]{$3$}
		\psfrag{g}[B][B]{$4$}
		\psfrag{h}[B][B]{$4$}
		\psfrag{j}[B][B]{$5$}
		\psfrag{k}[B][B]{$6$}
	\vspace*{8mm}
	\includegraphics[width=.95\linewidth]{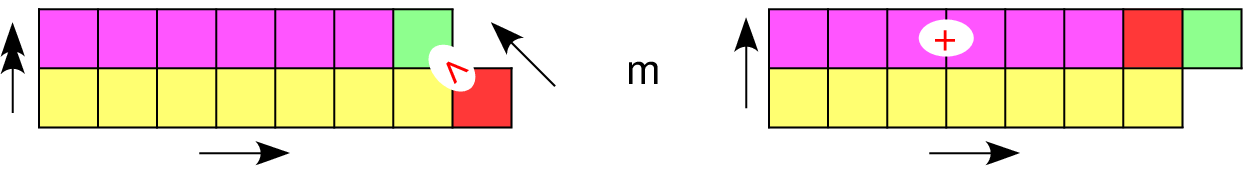}
	\vspace*{3mm}
	
	\includegraphics[width=.95\linewidth]{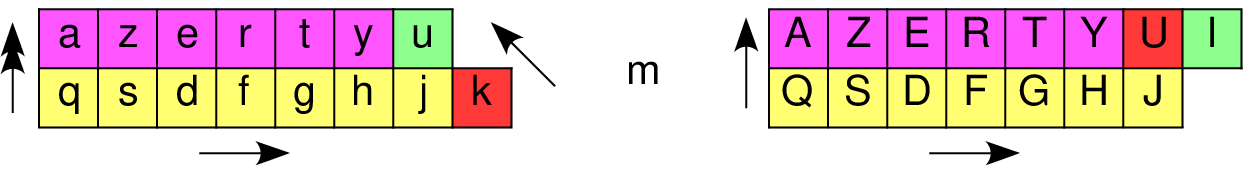}
	
		\psfrag{A}[B][B]{$1$}
		\psfrag{Z}[B][B]{$2$}
		\psfrag{E}[B][B]{$2$}
		\psfrag{R}[B][B]{$4$}
		\psfrag{T}[B][B]{$4$}
		\psfrag{Y}[B][B]{$6$}
		\psfrag{U}[B][B]{$6$}
		\psfrag{I}[B][B]{$6$}
		
		\psfrag{Q}[B][B]{$1$}
		\psfrag{S}[B][B]{$1$}
		\psfrag{D}[B][B]{$2$}
		\psfrag{F}[B][B]{$3$}
		\psfrag{G}[B][B]{$4$}
		\psfrag{H}[B][B]{$4$}
		\psfrag{J}[B][B]{$5$}
		
		\psfrag{a}[B][B]{$2$}
		\psfrag{z}[B][B]{$3$}
		\psfrag{e}[B][B]{$3$}
		\psfrag{r}[B][B]{$5$}
		\psfrag{t}[B][B]{$5$}
		\psfrag{y}[B][B]{$7$}
		\psfrag{u}[B][B]{$7$}
		
		\psfrag{q}[B][B]{$1$}
		\psfrag{s}[B][B]{$1$}
		\psfrag{d}[B][B]{$2$}
		\psfrag{f}[B][B]{$3$}
		\psfrag{g}[B][B]{$4$}
		\psfrag{h}[B][B]{$4$}
		\psfrag{j}[B][B]{$5$}
		\psfrag{k}[B][B]{$7$}
	\vspace*{8mm}
	\includegraphics[width=.95\linewidth]{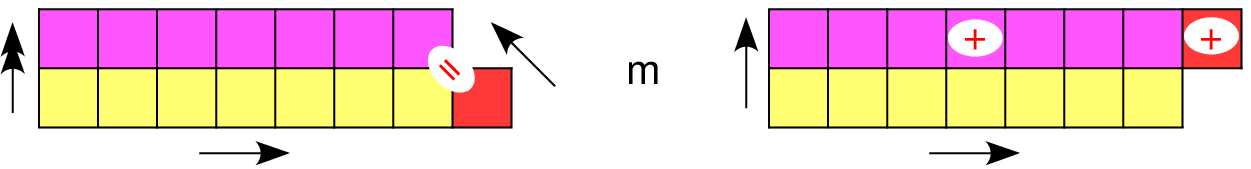}
	\vspace*{3mm}
	
	\includegraphics[width=.95\linewidth]{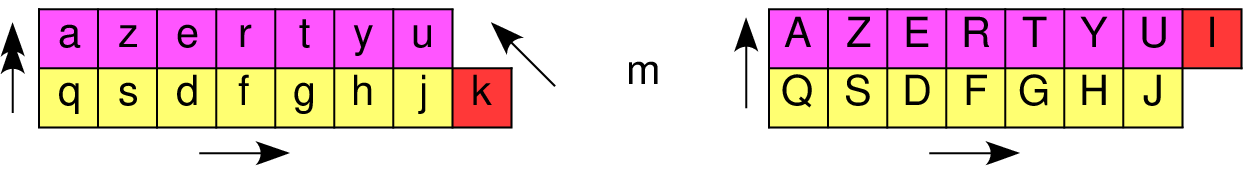}		
	\caption{The three cases of the bijection, from a Gog trapezoid to a Magog trapezoid. On the top line, $k=3$.}
	\label{gog_mag}
\end{figure}

\bigskip
\noindent\textbf{First case: $k\le n-2$.} We set
\begin{align*}
\col{violet}{m_{2,j}}&\de g_{2,j}-1\ \text{ for } 1\le j \le k-1\,; &&\col{rouge}{m_{2,k}}\de g_{1,k+1}\,;\qquad \col{vert}{m_{2,j}}\de g_{2,j-1}\ \text{ for } k+1 \le j \le n\,;\\
\col{jaune}{m_{1,j}}&\de g_{1,j}\ \text{ for } 1\le j \le k\,; &&\col{bleu}{m_{1,j}}\de g_{1,j+1}+2\ \text{ for }k+1\le j \le n-1\,.
\end{align*}

\medskip
\noindent\textbf{Second case: $k=n-1$ and $g_{1,n}<g_{2,n-1}$.} We set
\begin{align*}
\col{violet}{m_{2,j}}&\de g_{2,j}-1\ \text{ for } 1\le j \le n-2\,; &&\col{rouge}{m_{2,n-1}}\de g_{1,n}\,; &&\col{vert}{m_{2,n}}\de g_{2,n-1}\,;\\
\col{jaune}{m_{1,j}}&\de g_{1,j}\ \text{ for } 1\le j \le n-1\,.
\end{align*}

\medskip
\noindent\textbf{Third case: $k=n-1$ and $m_{1,n}=m_{2,n-1}$.} We set
\begin{align*}
\col{violet}{m_{2,j}}&\de g_{2,j}-1\ \text{ for } 1\le j \le n-1\, &&\col{rouge}{m_{2,n}}\de g_{1,n}-1\,;\\
\col{jaune}{m_{1,j}}&\de g_{1,j}\ \text{ for } 1\le j \le n-1\,.
\end{align*}

We now show that $\Psi_n(G)\in\M_n$. In the first and second case, the definition of~$k$ entails that $\col{violet}{m_{2,k-1}}=g_{2,k-1}-1 \le g_{1,k+1}=\col{rouge}{m_{2,k}}$, so that the horizontal inequalities hold. In the second case, we get the desired conclusion by noticing that $\col{rouge}{m_{2,n-1}}=g_{1,n}\le g_{2,n-1}-1\le n-1$ and $\col{vert}{m_{2,n}}=g_{2,n-1}\le n$. In the first case, by definition of~$k$, $\col{bleu}{m_{1,j}}=g_{1,j+1}+2 \le g_{2,j-1}=\col{vert}{m_{2,j}}$ for $k+1\le j \le n-1$ and, by vertical inequalities, $\col{vert}{m_{2,j}}=g_{2,j-1} \le j$ for $k+1\le j \le n$. Finally, still by definition of~$k$, $\col{rouge}{m_{2,k}}=g_{1,k+1}\le g_{1,k+2}\le g_{2,k}-2\le k-1$. This establishes the claim in the first case. The third case is straightforward.

\section{The previous mappings are inverses of each other}\label{secthm}

We now prove that the previous mappings are bijections.

\begin{thm}
The mappings $\Phi_n:\M_n\to\calG_n$ and $\Psi_n:\calG_n\to\M_n$ are bijections, which are inverse one from another.
\end{thm}

\begin{pre}
We have already established that $\Phi_n:\M_n\to\calG_n$ and $\Psi_n:\calG_n\to\M_n$. It remains to show that $\Psi_n\circ\Phi_n$ and $\Phi_n\circ\Psi_n$ are the identity on~$\M_n$ and~$\calG_n$, respectively. In fact, we will see that the three cases we distinguished are in correspondence via the bijection.

\medskip
\noindent\textbf{First case.} Let $M=(m_{i,j})\in\M_n$ be a Magog trapezoid that has a bug, and let~$k$ be its smallest bug. As in Section~\ref{secmg}, 
we define $(g_{ij})\de \Phi_n(M)$. We have $\col{violet}{g_{2,k-1}}=m_{2,k-1}+1\le m_{2,k}+1 = \col{rouge}{g_{1,k+1}}+1$ and, for $k+1\le j \le n-1$, $\col{bleu}{g_{1,j+1}}+1=m_{1,j}-1< m_{2,j}=\col{vert}{g_{2,j-1}}$ for $k+2\le j \le n-1$, so that
$$\max\big\{j\in\{2,\ldots,n-1\}\,: g_{2,j-1} \le g_{1,j+1}+1\big\}=k.$$
As the box moving procedure of Section~\ref{secgm} is clearly the inverse of that of Section~\ref{secmg}, we conclude that $(\Psi_n\circ\Phi_n)(M)=M$.

Let now $G=(g_{ij})\in\calG_n$ be such that the integer~$k$ defined by~\eqref{defk} is smaller than or equal to $n-2$. In order to conclude that $(\Phi_n\circ\Psi_n)(G)=G$, it is sufficient to show that~$k$ is the smallest bug of $(m_{i,j})\de\Psi_n(G)$. This is indeed the case as $\col{bleu}{m_{1,k+1}}=g_{1,k+2}+2> g_{1,k+1}+1=\col{rouge}{m_{2,k}}+1$ and, for $1\le j\le k-1$, $\col{jaune}{m_{1,j+1}}=g_{1,j+1}\le g_{2,j}=\col{violet}{m_{2,j}}+1$.

\medskip
\noindent\textbf{Second and third case.} Let $M=(m_{i,j})\in\M_n$ be a bug-free Magog trapezoid and $(g_{ij})\de \Phi_n(M)$. If we are in the second case, then $\col{violet}{g_{2,n-2}}=m_{2,n-2}+1\le m_{2,n-1}+1 = \col{rouge}{g_{1,n}}+1$, and, if we are in the third case, then $\col{violet}{g_{2,n-2}}=m_{2,n-2}+1\le m_{2,n}+1 = \col{rouge}{g_{1,n}}$, so that, in both cases,
$$\max\big\{j\in\{2,\ldots,n-1\}\,: g_{2,j-1} \le g_{1,j+1}+1\big\}=n-1.$$
We conclude as above that $(\Psi_n\circ\Phi_n)(M)=M$.

Let now $G=(g_{ij})\in\calG_n$ be such that the integer~$k$ defined by~\eqref{defk} is equal to $n-1$. We see that $(\Phi_n\circ\Psi_n)(G)=G$ by noticing that $(m_{i,j})\de\Psi_n(G)$ is bug-free as, for $1\le j\le n-1$, $\col{jaune}{m_{1,j+1}}=g_{1,j+1}\le g_{2,j}=\col{violet}{m_{2,j}}+1$.
\end{pre}

\section{Extension to \hr{$(\ell,n,2)$}{(l,n,2)} trapezoids and perspectives}\label{secln2}

\begin{figure}[ht]
		\psfrag{1}[][]{$4$}
		\psfrag{2}[][]{$5$}
		\psfrag{3}[][]{$6$}
		\psfrag{4}[][]{$7$}
		\psfrag{5}[][]{$8$}
		\psfrag{6}[][]{$9$}
		\psfrag{7}[][]{$10$}
		\psfrag{8}[][]{$11$}
		
		\psfrag{q}[B][B][.8]{\textcolor{red}{$g_{1,1}$}}
		\psfrag{s}[B][B][.8]{\textcolor{red}{$g_{1,2}$}}
		\psfrag{d}[B][B][.8]{\textcolor{red}{$g_{1,3}$}}
		\psfrag{f}[B][B][.8]{\textcolor{red}{$g_{1,4}$}}
		\psfrag{g}[B][B][.8]{\textcolor{red}{$g_{1,5}$}}
		\psfrag{h}[B][B][.8]{\textcolor{red}{$g_{1,6}$}}
		\psfrag{j}[B][B][.8]{\textcolor{red}{$g_{1,7}$}}
		\psfrag{k}[B][B][.8]{\textcolor{red}{$g_{1,8}$}}
		
		\psfrag{a}[B][B][.8]{\textcolor{red}{$g_{2,1}$}}
		\psfrag{z}[B][B][.8]{\textcolor{red}{$g_{2,2}$}}
		\psfrag{e}[B][B][.8]{\textcolor{red}{$g_{2,3}$}}
		\psfrag{r}[B][B][.8]{\textcolor{red}{$g_{2,4}$}}
		\psfrag{t}[B][B][.8]{\textcolor{red}{$g_{2,5}$}}
		\psfrag{y}[B][B][.8]{\textcolor{red}{$g_{2,6}$}}
		\psfrag{u}[B][B][.8]{\textcolor{red}{$g_{2,7}$}}
		
		\psfrag{Q}[B][B][.8]{\textcolor{red}{$m_{1,1}$}}
		\psfrag{S}[B][B][.8]{\textcolor{red}{$m_{1,2}$}}
		\psfrag{D}[B][B][.8]{\textcolor{red}{$m_{1,3}$}}
		\psfrag{F}[B][B][.8]{\textcolor{red}{$m_{1,4}$}}
		\psfrag{G}[B][B][.8]{\textcolor{red}{$m_{1,5}$}}
		\psfrag{H}[B][B][.8]{\textcolor{red}{$m_{1,6}$}}
		\psfrag{J}[B][B][.8]{\textcolor{red}{$m_{1,7}$}}
		
		\psfrag{A}[B][B][.8]{\textcolor{red}{$m_{2,1}$}}
		\psfrag{Z}[B][B][.8]{\textcolor{red}{$m_{2,2}$}}
		\psfrag{E}[B][B][.8]{\textcolor{red}{$m_{2,3}$}}
		\psfrag{R}[B][B][.8]{\textcolor{red}{$m_{2,4}$}}
		\psfrag{T}[B][B][.8]{\textcolor{red}{$m_{2,5}$}}
		\psfrag{Y}[B][B][.8]{\textcolor{red}{$m_{2,6}$}}
		\psfrag{U}[B][B][.8]{\textcolor{red}{$m_{2,7}$}}
		\psfrag{I}[B][B][.8]{\textcolor{red}{$m_{2,8}$}}		
		
		\psfrag{n}[][]{$(3,8,2)$-Gog trapezoid}
		\psfrag{m}[][]{$(3,8,2)$-Magog trapezoid}
	\centering\includegraphics[width=.95\linewidth]{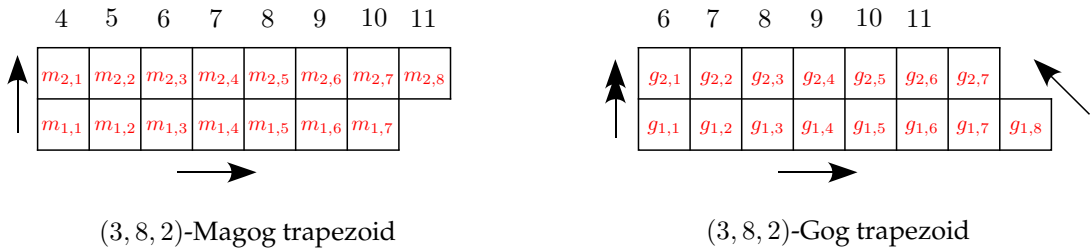}
	\caption{Definition of $(\ell,n,2)$ trapezoids.}
	\label{def2}
\end{figure}

Our bijection can trivially be extended to $(\ell,n,2)$ trapezoids, where $\ell\ge 0$ is an integer. Here, an \emph{$(\ell,n,2)$-Magog trapezoid} is defined as an $(n,2)$-Magog trapezoid, with the difference that item~\ref{defmii} of Definition~\ref{defm} is replaced by
\begin{enumerate}[label=(\textit{ii'})]
	\item $m_{1,j}\le m_{2,j}\le j+\ell$ for all $j\in\{1,\ldots, n-1\}$ and $m_{2,n}\le n+\ell$\,.
\end{enumerate}
See Figure~\ref{def2}. Similarly, an \emph{$(\ell,n,2)$-Gog trapezoid} is defined as an $(n,2)$-Gog trapezoid with the difference that item~\ref{defgii} of Definition~\ref{defg} is replaced by
\begin{enumerate}[label=(\textit{ii'})]
	\item $g_{1,j}< g_{2,j}<j+2+\ell$ for all $j\in\{1,\ldots, n-1\}$\,;
\end{enumerate}

For any $\ell\ge 1$ and $n\ge 3$, the mappings~$\Phi_n$ and~$\Psi_n$ can be extended without any differences in the construction into bijections between the set of $(\ell,n,2)$-Magog trapezoids and the set of $(\ell,n,2)$-Gog trapezoids. The proofs can be copied almost verbatim, the only thing to do is to add~$\ell$ whenever we use one of the bounds changed by these definitions.

\bigskip

Unfortunately, as of today, we did not manage to extend this bijection to $(n,3)$ trapezoids. The mapping~$\Phi_n$ exchanges the sizes of two consecutive rows so that one could think that, in the case of $(n,3)$ trapezoids, we would need to apply a similar operation several times in order to pass from a Magog to a Gog trapezoid. Unfortunately, whenever a third row is present, we cannot slide the boxes of two consecutive rows without breaking the rules. This question remains under investigation.

\bigskip

It has also been brought to our attention that our construction bears some intriguing similarities with a construction used by Krattenthaler \cite[Section~2]{krattenthaler89} in order to show the $q$-log-concavity of Gaussian binomial coefficients. In the latter construction, two rows of strictly increasing integers are considered and one carefully chosen entry of the second row is removed from it and inserted in the first row after addition of a constant. A major difference between both constructions lies in the fact that the entries in the latter one only satisfy monotonicity relations in one direction (along rows) so that the objects are less constrained.